\documentclass{elsart3-1}

\usepackage{ae}
\usepackage{aecompl}
\usepackage[cm]{aeguill}
\usepackage{amsmath, amssymb, amsfonts}
\usepackage[english,francais]{babel}
\usepackage{dsfont}

\newtheorem{theorem}{Theorem}[section]

\newtheorem{e-proposition}[theorem]{Proposition}

\newtheorem{e-definition}[theorem]{Definition\rm}

\renewcommand{\theequation}{\arabic{equation}}
\def\og{\leavevmode\raise.3ex\hbox{$\scriptscriptstyle\langle\!\langle$~}}
\def\fg{\leavevmode\raise.3ex\hbox{~$\!\scriptscriptstyle\,\rangle\!\rangle$}}

\def\E{\mathbb{E}}

\def\1{\mathds{1}}

\newenvironment{disarray}
{\everymath{\displaystyle\everymath{}}\array}
{\endarray}

\begin{document}
\centerline{Statistics}
\begin{frontmatter}


\selectlanguage{english}
\title{Rates of convergence for nonparametric deconvolution}


\selectlanguage{english}
\author[authorlabel1]{Claire LACOUR},
\ead{lacour@math-info.univ-paris5.fr}

\address[authorlabel1]{Laboratoire MAP 5,
Universit\'e Paris 5, 
45 rue des Saints-P\`eres,
75270 Paris Cedex 06}


\medskip

\begin{abstract}
\selectlanguage{english}
This Note presents original rates of convergence for the deconvolution problem. We assume
that both the estimated density and noise density are supersmooth and we compute the 
risk for two kinds of estimators.

\vskip 0.5\baselineskip

\selectlanguage{francais}
\noindent{\bf R\'esum\'e} \vskip 0.5\baselineskip \noindent
{\bf Vitesses de convergence en d\'econvolution nonparam\'etrique.}
Cette Note pr\'esente des vitesses de convergence originales pour le probl\`eme de
d\'econvolution. On suppose que la densit\'e estim\'ee ainsi que la densit\'e du bruit 
sont \og supersmooth\fg ~ et on calcule le risque pour deux types d'estimateurs.

\end{abstract}
\end{frontmatter}

\selectlanguage{francais}
\section*{Version fran\c{c}aise abr\'eg\'ee}


On consid\`ere le probl\`eme de d\'econvolution suivant : 
\begin{equation*}
Y_i=X_i+\varepsilon_i  \qquad \qquad i=1,\dots, n
\end{equation*}
o\`u les $X_i$ sont des variables al\'eatoires ind\'ependantes et identiquement distribu\'ees de densit\'e inconnue $g$ 
et les variables al\'eatoires $\varepsilon_i$ sont i.i.d de densit\'e connue $f_\varepsilon$.
Les suites 
$(X_i)$ et  $(\varepsilon_i)$ sont de plus suppos\'ees ind\'ependantes.
L'objectif est d'estimer $g$ \`a partir des donn\'ees $Y_1,\dots, Y_n$.

Le cadre d'hypoth\`eses est le suivant. 
Notons, pour toute fonction $u$, $u^*$  
la transform\'ee de Fourier de $u$: $u^*(x)=\int e^{ixt}u(t)dt$.
On suppose que le bruit est tel que pour tout $x$ de $\mathbb{R}$, $f_\varepsilon^*(x)\neq 0$ 
et qu'il satisfait l'hypoth\`ese suivante :
\begin{description}
\item[N1.] Il existe $s\geq0, b\geq 0, \gamma \in\mathbb{R}$ ($\gamma>0$ si $s=0$)
et $k_0, k_1>0$ tels que
$$k_0(x^2+1)^{-\gamma/2}\exp(-b|x|^s)\leq |f_\varepsilon^*(x)|\leq k_1(x^2+1)^{-\gamma/2}\exp(-b|x|^s)$$
\end{description}
On suppose de plus que $g$ appartient \`a l'espace 
$$\mathcal{A}_{\delta,r,a}(L)=\{g \text{ densit\'e sur } \mathbb{R} \text{ et } 
\int |g^*(x)|^2(x^2+1)^{\delta}\exp(2a|x|^r)\leq L\}$$
avec $r\geq0, a\geq 0, \delta \in\mathbb{R}$ ($\delta>1/2$ si $r=0$), $L>0$. Lorsque $r$ est
strictement positif, la fonction est dite superr\'eguli\`ere et  ordinairement r\'eguli\`ere sinon, 
la terminologie \'etant la m\^eme pour le bruit.

Ce probl\`eme a \'et\'e largement \'etudi\'e pour une fonction $g$ appartenant \`a
un espace de Sobolev ou de H\"{o}lder (i.e. $r=0$) : on peut citer entre autres \cite{5}, \cite{6}, \cite{7,8}, 
\cite{9}, \cite{15}. Les vitesse de convergences m\'ediocres (en puissances de $\ln n$) obtenues lorsque le bruit 
est  superr\'egulier (et en particulier pour les distributions gaussiennes) ont conduit \`a consid\'erer 
des fonctions $g$ \'egalement superr\'eguli\`eres. 
En premier \cite{12} et plus r\'ecemment \cite{butucea1}, \cite{butucea} et \cite{comte} 
ont \'etudi\'e des estimateurs dans ce contexte.


Dans cette Note, nous fournissons des vitesses de convergence exactes et explicites, 
m\^eme dans le cas $r>0, s>0$ pour lequel jusqu'\`a maintenant les vitesses n'\'etaient 
donn\'ees que de fa\c{c}on  implicite, except\'e dans des cas tr\`es particuliers. 
Ces vitesses sont calcul\'ees pour deux types d'estimateurs.


L'estimateur \`a noyau classique est le suivant : 
$$\hat{g}_n(x)=\frac{1}{nh} \sum_{i=1}^n K\left(\frac{x-Y_i}{h}\right) $$
 o\`u $K$ est d\'efini comme la transform\'ee de Fourier inverse de
$K^*(x)=\1_{\{|x|\leq 1\}}/f_\varepsilon^*(x/h)$. 
Le risque ponctuel (not\'e MSE) et le risque int\'egr\'e (not\'ee MISE) sont \'etablis dans 
\cite{butucea} (voir la proposition \ref{prop1} ci-dessous).
Cet estimateur a l'avantage d'\^etre optimal au sens minimax asymptotique fin (voir \cite{butucea})
mais ne fournit d'estimateur adaptatif que dans des cas particuliers. C'est pourquoi nous nous
int\'eressons \'egalement \`a l'estimateur de projection introduit dans \cite{comte}.

Soit $\varphi(x)=\sin(\pi x)/(\pi x)$ et $\varphi_{m,j}(x)=\sqrt{L_m}\varphi(L_mx-j)$.
En posant $u_t(x)=(1/2\pi)(t^*(-x)/f_\varepsilon^*(x))$,  l'estimateur de projection est d\'efini par
$$ \hat{g}_m(x)=\sum_{|j|\leq K_n}\hat{a}_{m,j}\varphi_{m,j}\quad\quad \text{   o\`u  } 
\hat{a}_{m,j}=\frac{1}{n} \sum_{i=1}^n u^*_{\varphi_{m,j}}(Y_i)$$
Pour cet estimateur, le calcul du rique int\'egr\'e est pr\'esent\'e dans \cite{comte} (voir la proposition 
\ref{prop2} ci-dessous).
Un estimateur adaptatif peut ensuite \^etre d\'efini en utilisant une technique de s\'election de mod\`eles 
(voir \cite{comte} pour les d\'etails).


On peut observer que les deux estimateurs ont les m\^emes vitesses de convergence, qui
s'obtiennent en minimisant l'ordre du risque. 
On est alors ramen\'e \`a r\'esoudre l'\'equation suivante
\begin{equation*}
  \exp(\frac{2b}{h^s}+\frac{2a}{h^r})h^\alpha =O(n) 
\end{equation*}
o\`u $\alpha=r-2\delta-2\gamma-1$ lorsque l'on consid\`ere l'erreur int\'egr\'ee et
$\alpha=-2\delta-2\gamma+(s-1)_+$ lorsque l'on consid\`ere l'erreur ponctuelle. Dans la plupart des cas,
la solution de cette \'equation est bien connue (voir Tables \ref{rates} et \ref{rates2}).
Seul le cas $r>0$, $s>0$ n'a pas \'et\'e compl\`etement r\'esolu. C'est ce cas qui est \'etudi\'e ici
: les vitesses sont donn\'ees explicitement dans le th\'eor\`eme \ref{main}.

Trois cas sont \`a diff\'erencier : $r=s$, $ r<s$ et $ r>s$.
Si $r$ est \'egal \`a $s$, les vitesses (ponctuelles et int\'egr\'ees) sont en 
$n^{-a/(a+b)}$ modifi\'ees par un facteur logarithmique.
Si $r$ est diff\'erent de $s$, on observe un ph\'enom\`ene original :
les vitesses d\'ependent de l'intervalle $]k/(k+1),(k+1)/(k+2)], k\in \mathbb{N}$ auquel appartient $r/s$ ou $s/r$.
Si $r$ est strictement inf\'erieur \`a $s$ (biais dominant), le terme principal est d'ordre 
$\exp[b_0(\ln n)^{r/s}]$ avec $b_0=-2a/(2b)^{r/s}$, et si $r$ est strictement sup\'erieur \`a $s$ 
(variance dominante), le terme principal est d'ordre 
$\exp[-d_0(\ln n)^{s/r}]/n$ avec $d_0=-2b/(2a)^{s/r}$.

Ainsi ces vitesses originales d\'ecroissent plus vite que n'importe quelle 
fonction logarithmique. De plus elles sont optimales lorque les bornes inf\'erieures correspondantes sont connues,
c'est-\`a-dire $r=s=1$ (voir \cite{tsy}) et $r<s$ (voir\cite{butucea}).
Il est \`a noter qu'\'etant donn\'e la complexit\'e de ces vitesses, il 
est r\'eellement int\'eressant de construire des estimateurs adaptatifs, c'est-\`a-dire
des estimateurs dont le risque atteint automatiquement les vitesses minimax.

\selectlanguage{english}
\section{Introduction}

We consider the following deconvolution problem:
\begin{equation*}
Y_i=X_i+\varepsilon_i  \qquad \qquad i=1,\dots, n
\end{equation*}
where the $X_i$'s are independent and identically distributed random variables with an unknown density $g$ and
the random variables $\varepsilon_i$ are i.i.d with known density $f_\varepsilon$. Moreover
$(X_i)$ and  $(\varepsilon_i)$ are independent.
The aim is to estimate $g$ from data $Y_1,\dots, Y_n$.

The hypothesis framework is the following.
Denote, for all function $u$, $u^*$ the Fourier transform of $u$: $u^*(x)=\int e^{ixt}u(t)dt$.
We suppose that noise is such that for all $x$ in $\mathbb{R}$, $f_\varepsilon^*(x)\neq 0$ 
and that it satisfies the following assumption:
\begin{description}
\item[N1.] There exist $s\geq0, b\geq 0, \gamma \in\mathbb{R}$ ($\gamma>0$ if $s=0$)
and $k_0, k_1>0$ such that
$$k_0(x^2+1)^{-\gamma/2}\exp(-b|x|^s)\leq |f_\varepsilon^*(x)|\leq k_1(x^2+1)^{-\gamma/2}\exp(-b|x|^s)$$
\end{description}
We assume that $g$ belongs to the space 
$$\mathcal{A}_{\delta,r,a}(L)=\{g \text{ is a probability density on } \mathbb{R} \text{ and } 
\int |g^*(x)|^2(x^2+1)^{\delta}\exp(2a|x|^r)\leq L\}$$
with $r\geq0, a\geq 0, \delta \in\mathbb{R}$ ($\delta>1/2$ if $r=0$), $L>0$. When $r>0$
the function is known as supersmooth, and as ordinary smooth else. The terminology is the same
for noise.

This problem has been extensively studied for a function $g$ belonging to a Sobolev 
or H\"{o}lder class (i.e. $r=0$): see among others \cite{5}, \cite{6}, \cite{7,8}, 
\cite{9}, \cite{15}. The bad rates of convergence (power of $\ln n$) for 
supersmooth noise (and then in particular for Gaussian distributions) lead to consider supersmooth 
functions. First \cite{12} and more recently \cite{butucea1}, \cite{butucea} and \cite{comte} 
studied estimators in this context.

The contribution of this paper is to provide exact and explicit rates of convergence, 
even in the case $r>0$ and $s>0$ where up to now the rates were not explicitly available except in
 very particular cases.

\section{Estimators and preliminar results}

\subsection{Estimators}
The classical kernel estimator is the following: 
$$\hat{g}_n(x)=\frac{1}{nh} \sum_{i=1}^n K\left(\frac{x-Y_i}{h}\right) $$
 where $K$ is the function defined as the inverse Fourier transform of 
$K^*(x)=\1_{\{|x|\leq 1\}}/f_\varepsilon^*(x/h)$. 
The pointwise mean squared error (denoted by MSE) and mean integrated squared error (denoted by MISE)
are established in \cite{butucea}:
\begin{e-proposition} \label{prop1} If $g$ belongs to $\mathcal{A}_{\delta,r,a}(L)$, then under Assumption N1,\\
$MISE=\E\|g-\hat{g}_n\|_2^2=O\left(h^{2\delta}\exp(-2a/h^r)+\cfrac{h^{s-1-2\gamma}}{n}\exp(2b/h^s)\right)\quad$  and  \\
$MSE=\E|g(x)-\hat{g}_n(x)|^2=O\left(h^{2\delta+r-1}\exp(-2a/h^r)+\min(1,h^{s-1})\cfrac{h^{s-1-2\gamma}}{n}\exp(2b/h^s)\right)$
\end{e-proposition}
This estimator has the advantage to be optimal in the sharp asymptotic minimax sense (see \cite{butucea})
but provides an adaptive estimator only in particular cases. That is why we present the projection 
estimator introduced in \cite{comte}.

Let $\varphi(x)=\sin(\pi x)/(\pi x)$ and $\varphi_{m,j}(x)=\sqrt{L_m}\varphi(L_mx-j)$.
Consider $u_t(x)=(1/2\pi)(t^*(-x)/f_\varepsilon^*(x)).$ Then the projection estimator is defined by
$$ \hat{g}_m(x)=\sum_{|j|\leq K_n}\hat{a}_{m,j}\varphi_{m,j}\quad\quad \text{   where  } 
\hat{a}_{m,j}=\frac{1}{n} \sum_{i=1}^n u^*_{\varphi_{m,j}}(Y_i)$$
For this estimator, the following result is proved in \cite{comte}:
\begin{e-proposition} \label{prop2} Assume that $f_\varepsilon\in L^2$ (i.e. $\gamma>1/2$ when $s=0$)
and that $g$ is a $L^2$ function which verifies $\int x^2g^2(x)dx\leq M$.
If $g$ belongs to $\mathcal{A}_{\delta,r,a}(L)$, then under Assumption N1,\\
$MISE=\E\|g-\hat{g}_m\|_2^2=O\left(L_m^{-2\delta}\exp(-2a\pi^rL_m^r)+\cfrac{L_m^{2\gamma+1-s}}{n}\exp(2b\pi^sL_m^s)\right)$
\end{e-proposition}
Then, an adaptive estimator can be defined using  a model selection method (see \cite{comte} for details).

\subsection{Rates of convergence}

We can observe that both estimators have the same $L^2$ rate of convergence (take 
$h^{-1}=\pi L_m$). To compute this rate, we have to minimize the 
risk orders in $h$ (or $L_m$). By setting to zero the derivative of this quantity we obtain the equation 
\begin{equation}
  \exp(\frac{2b}{h^s}+\frac{2a}{h^r})h^\alpha =O(n) \label{equa}
\end{equation}
where $\alpha=r-2\delta-2\gamma-1$ if we consider the integrated error and 
$\alpha=-2\delta-2\gamma+(s-1)_+$ if we consider the pointwise error. In most cases,
the solution of this equation is well known and leads to the following tables where 
different regularities for $g$ and $f_\varepsilon$ are examined:

\begin{table}[!h]
\begin{minipage}[t]{.05\linewidth}
\quad
\end{minipage}
\begin{minipage}[t]{.3\linewidth}
 $$\begin{disarray}{ccc}
\cline{2-3}
& \quad s=0 \quad & \quad s>0 \quad\\
\hline
\quad r=0 \quad & \quad n^{-\frac{2\delta}{2\delta+2\gamma+1}} \quad & 
\quad (\ln n)^{-\frac{2\delta}{s}} \quad \\
\hline
\quad r>0 \quad &\quad \frac{(\ln n)^{\frac{2\gamma+1}{r}}}{n} \quad & \quad \text{Theorem \ref{main}}\quad  \\
\hline
\end{disarray}$$
\caption{Rates of convergence for the MISE.}\label{rates}
\end{minipage} 
\begin{minipage}[t]{.15\linewidth}
\quad
\end{minipage}
\begin{minipage}[t]{.3\linewidth}
 $$\begin{disarray}{ccc}
\cline{2-3}
 &\quad s=0\quad &\quad s>0\quad \\
\hline
\quad r=0 \quad &\quad n^{\frac{1-2\delta}{2\delta+2\gamma}} \quad & \quad (\ln n)^{\frac{1-2\delta}{s}}\quad\\
\hline
\quad r>0 \quad & \quad \frac{(\ln n)^{\frac{2\gamma+1}{r}}}{n}\quad & \quad \text{Theorem \ref{main}}\quad \\
\hline
\end{disarray}$$
\caption{Rates of convergence for the MSE.}\label{rates2}
\end{minipage}
\end{table}

Except for the bottom right cells (to be completed in the next section), these rates are known 
to be optimal minimax rates: see \cite{7} and \cite{butucea1} for the lower bounds.

\section{Results}

The rates of convergence in the case $(r>0, s>0)$ depend on the integer $k$ such that 
$r/s$ or $s/r$ belongs to the interval $(k/(k+1),(k+1)/(k+2)]$:
\begin{theorem}\label{main}
We assume $r>0$ and $s>0$. Let $k\in \mathbb{N}$ and $\lambda=\mu^{-1}=r/s$. Then
\begin{description}
\item $\bullet$ if $r=s$, if $\xi=[2\delta b+(s-2\gamma-1)a]/[(a+b)s]$
$$\displaystyle MISE=O\left(n^{-{a}/(a+b)}(\ln n)^{-\xi}\right);$$
$$MSE=O\left(n^{-{a}/(a+b)}(\ln n)^{-\xi+\frac{(1-s)_+b}{(a+b)s}}\right)$$
\item $\bullet$ if $r<s$ and $\cfrac{k}{k+1}<\lambda\leq\cfrac{k+1}{k+2}$, there exist reals $b_i$
such that  
$$MISE=O\left((\ln n)^{-{2\delta}/{s}}\exp[\sum_{i=0}^k b_i(\ln n)^{(i+1)\lambda-i}]\right);$$
$$MSE=O\left((\ln n)^{(-2\delta-r+1)/{s}}\exp[\sum_{i=0}^k b_i(\ln n)^{(i+1)\lambda-i}]\right)$$
\item $\bullet$ if $r>s$ and $\cfrac{k}{k+1}<\mu\leq\cfrac{k+1}{k+2}$, there exist reals $d_i$
such that,  
$$MISE=O\left(\frac{(\ln n)^{({1+2\gamma-s})/{r}}}{n}\exp[-\sum_{i=0}^k d_i(\ln n)^{(i+1)\mu-i}]\right);$$
$$MSE=O\left(\frac{(\ln n)^{({1+2\gamma-s-(s-1)_+})/{r}}}{n}
\exp[-\sum_{i=0}^k d_i(\ln n)^{(i+1)\mu-i}]\right)$$
\end{description}
\end{theorem}

The coefficients $b_i$ and $d_i$ are computable, see Section \ref{proof} for the exact form of reals $b_i$.
Notice that these original rates have the property to decrease faster than any logarithmic function. 
Moreover, they are optimal in the cases where the corresponding lower bounds are known, i.e $r=s=1$ 
(see \cite{tsy}) and $r<s$ (see\cite{butucea}). We can also remark that, 
given the complexity of these rates, it is woth finding adaptive estimators, i.e.
estimators whose risk automatically achieves the minimax rates.

\section{Proof of Theorem \ref{main}}\label{proof}

\begin{description}
\item $\bullet$ If $r=s$, we check that $h^*=(2a+2b)^{1/s}\left(\ln n+\frac{\alpha}{s}\ln\ln n\right)^{-1/s}$
satisfies Equation \eqref{equa}. The corresponding risks are easily obtained.
\item $\bullet$ If $r<s$ and $\cfrac{k}{k+1}<\lambda\leq\cfrac{k+1}{k+2}$, let
$h*=(2b)^{1/s}[\ln n+\frac{\alpha}{s}\ln\ln n+\sum_{i=0}^k b_i(\ln n)^{(i+1)\lambda-i}]^{-1/s}$
$$\begin{disarray}{rcl}
 \exp[\frac{2B}{h^{*s}}+\frac{2A}{h^{*r}}+\alpha \ln h^*]&=&Kn\exp[
 \frac{2a}{(2b)^\lambda}(\ln n)^\lambda(1+\lambda u_n+..+\frac{\lambda(\lambda-1)..(\lambda-k)}
 {(k+1)!}u_n^{k+1}+o(u_n^{k+1}))\\
& & +\sum_{i=0}^k b_i(\ln n)^{(i+1)\lambda-i}-\frac{\alpha}{s}\ln(1+o(1))]   
\end{disarray}$$
with $u_n:=\frac{\alpha}{s}\frac{\ln\ln n}{\ln n}
   +\sum_{i=0}^k b_i(\ln n)^{(i+1)\lambda-(i+1)}$\\
By noting that for $2\leq j \leq k+1$,
$u_n^j=\sum_{i=j}^{k+1}\sum_{p_1+..+p_j=i}b_{p_1-1}..b_{p_j-1}(\ln n)^{(\lambda-1)i}+o((\ln n)^{(\lambda-1)(k+1)})$,
we compute 
$$\exp[\frac{2b}{h^{*s}}+\frac{2a}{h^{*r}}+\alpha \ln h^*]=Kn \exp[o(1)+\sum_{i=0}^{k+1}M_i(\ln n)^{(i+1)\lambda-i}+
  o((\ln n)^{(k+2)\lambda-(k+1)})]$$
with $\displaystyle M_0=b_0+\cfrac{2a}{(2b)^\lambda}$ and for $1\leq i \leq k+1$,
$\displaystyle M_i=b_i+\cfrac{2a}{(2b)^\lambda}\sum_{j=1}^i\cfrac{\lambda\dots(\lambda-j+1)}{j!}
\sum_{p_1+..+p_j=i}b_{p_1-1}..b_{p_j-1}$\\
by denoting $b_{k+1}=0$.
The $b_i$'s are recursively defined by $M_0=...=M_k=0$ and thus Equation \eqref{equa} is verified.
The risks are then computed by using the same tools.
\item $\bullet$ The reasoning is the same in the case $r>s$

\end{description}

\end{document}